\documentclass[12pt]{article}
\usepackage{amssymb}

\newcommand{\text}{\mbox}

\newtheorem{theorem}{Theorem}[section]
\newtheorem{proposition}{Proposition}[section]
\newtheorem{lemma}{Lemma}[section]
\newtheorem{corollary}{Corollary}[section]
\newtheorem{definition}{Definition}[section]

\begin{document}
\author{Frank Oertel\\Zurich Financial Services, Reinsurance\\
Corporate Development and Risk Management\\
General Guisan-Quai 26\\
8022 Zurich, Switzerland\\
\textit{E-mail address: oertel@e-math.ams.org}}
%\thanks{E-mail address: frank.oertel@zurichre.com}
\title{ The principle of local reflexivity for operator ideals and
its implications}
\maketitle

\begin{abstract}
\noindent We present a survey of past research activities and current
results in constructing a mathematical framework describing the
principle of local reflexivity for operator ideals and reveal
further applications involving operator ideal products consisting
of operators which factor through a Hilbert space.
\end{abstract}

\section{Introduction}

This survey paper is devoted to the presentation of the author`s research
results related to the investigation of operator ideals $(\frak{A},\mathbf{A}%
)$ which allow a transfer of the norm estimation in the classical principle
of local reflexivity to their ideal (quasi--)norm $\mathbf{A}$. This
research originated from the objective to facilitate the search for a
non--accessible maximal normed Banach ideal (which is the same as a
non--accessible finitely generated tensor norm in the sense of Grothendieck)
and lead to the dissertation \cite{oe1} in 1990. Later, in 1993, Pisier
constructed such a counterexample (cf. \cite[31.6.]{df}). Since each
right--accessible maximal Banach ideal $(\frak{A},\mathbf{A})$ even
satisfies such a \textit{principle of local reflexivity for operator ideals }%
(called $\frak{A}$--\textit{LRP}), Pisier`s counterexample of a
non--accessible maximal Banach ideal naturally lead to the search
for counterexamples of maximal Banach ideals
$({\frak{A}}_0,{\mathbf{A}}_0)$ which even do not satisfy the
${\frak{A}}_0$-LRP, implying surprising relations between the
existence of
a norm on product operator ideals of type $\frak{B}\circ \frak{L}_2$ (where $%
\frak{L}_2$ denotes the class of all operators which factor through a
Hilbert space), the extension of finite rank operators with respect to a
suitable operator ideal norm and the principle of local reflexivity for
operator ideals (cf. \cite{oe5}). The basic objects, connecting these
different aspects, are product operator ideals with property (I) and
property (S), introduced by Jarchow and Ott (see \cite{jo}). In the widest
sense, a product operator ideal $\frak{A}\circ \frak{B}$ has the property
(I), if
\[
(\frak{A}\circ \frak{B})\cap \frak{F}=(\frak{A}\cap \frak{F})\circ \frak{B}
\]
and the property (S), if
\[
(\frak{A}\circ \frak{B})\cap \frak{F}=\frak{A}\circ (\frak{B}\cap \frak{F})% \text{,}
\]
(where $\frak{F}$ denotes the class of all finite rank operators) so that
each finite rank operator in $\frak{A}\circ \frak{B}$ is the composition of
two operators, one of which is of finite rank. Since each operator ideal
which contains $\frak{L}_2$ as a factor, has both, the property (I) and the
property (S), Hilbert space factorization crystallized out as a fundamental
key in these investigations.

\section{The framework}

In this section, we introduce the basic notation and terminology which we
will use throughout in this paper. We only deal with Banach spaces and most
of our notations and definitions concerning Banach spaces and operator
ideals are standard. We refer the reader to the monographs \cite{df}, \cite
{djt} and \cite{p1} for the necessary background in operator ideal theory
and the related terminology. Infinite dimensional Banach spaces over the
field $\Bbb{K\in \{R},\Bbb{C\}}$ are denoted throughout by $W,X,Y$ and $Z$
in contrast to the letters $E,F$ and $G$ which are used for finite
dimensional Banach spaces only. The space of all operators (continuous
linear maps) from $X$ to $Y$ is denoted by $\frak{L}$(X,Y), and for the
identity operator on $X$, we write $Id_X$. The collection of all finite rank
(resp. approximable) operators from $X$ to $Y$ is denoted by $\frak{F}$(X,Y)
(resp. $\overline{\frak{F}}$(X,Y)), and $\frak{E}$(X,Y) indicates the
collection of all operators, acting between finite dimensional Banach spaces
$X$ and $Y$ (elementary operators). The dual of a Banach space $X$ is
denoted by $X^{\prime }$, and $X^{\prime \prime }$ denotes its bidual $%
(X^{\prime })^{\prime }$. If $T\in \frak{L}$(X,Y) is an operator, we
indicate that it is a metric injection by writing $T:X\stackrel{1}{%
\hookrightarrow }Y$, and if it is a metric surjection, we write $T:X%
\stackrel{1}{\twoheadrightarrow }Y$. If $X$ is a Banach space, $E$ a finite
dimensional subspace of $X$ and $K$ a finite codimensional subspace of $X$,
then $B_X:=\{x\in X:\Vert x\Vert \leq 1\}$ denotes the closed unit ball, $%
J_E^X$ $:E\stackrel{1}{\hookrightarrow }X$ the canonical metric injection
and $Q_K^X:X\stackrel{1}{\twoheadrightarrow }X\diagup K$ the canonical
metric surjection. Finally, $T^{\prime }\in {\frak{L}}(Y^{\prime },X^{\prime
}) $ denotes the dual operator of $T\in \frak{L}$(X,Y).

If $({\frak{A}},\mathbf{A})$ and $({\frak{B}},\mathbf{B})$ are given
quasi--Banach ideals, we will use throughout the shorter notation $({\frak{A}}%
^d,\mathbf{A}^d)$ for the dual ideal and the abbreviation ${\frak{A}}\stackrel{%
1}{=}{\frak{B}}$ for the isometric equality $(\frak{A},\mathbf{A})=(\frak{B},%
\mathbf{B})$. We write $\frak{A}\subseteq \frak{B}$ if, regardless of the
Banach spaces $X$ and $Y$, we have ${\frak{A}}(X,Y)\subseteq {\frak{B}}(X,Y)$.
If $X_0$ is a fixed Banach space, we write ${\frak{A}}(X_0,\cdot )\subseteq
{\frak{B}}(X_0,\cdot )$ (resp. ${\frak{A}}(\cdot ,X_0)\subseteq {\frak{B}}(\cdot
,X_0)$) if, regardless of the Banach space $Z$ we have ${\frak{A}}%
(X_0,Z)\subseteq {\frak{B}}(X_0,Z)$ (resp. ${\frak{A}}(Z,X_0)\subseteq {\frak{B}}%
(Z,X_0)$). The metric inclusion $(\frak{A},\mathbf{A})\subseteq (\frak{B},%
\mathbf{B})$ is often shortened by $\frak{A}\stackrel{1}{\subseteq }\frak{B}$%
. If ${\mathbf{B}}(T)\leq {\mathbf{A}}(T)$ for all finite rank (resp.
elementary) operators $T\in \frak{F}$ (resp. $T\in \frak{E}$), we sometimes
use the abbreviation $\frak{A}\stackrel{\frak{F}}{\subseteq }\frak{B}$
(resp. $\frak{A}\stackrel{\frak{E}}{\subseteq }\frak{B}$).

First we recall the basic notions of Grothendieck's metric theory of tensor
products (cf., eg., \cite{df}, \cite{gl}, \cite{gr}, \cite{l}), which
together with Pietsch's theory of operator ideals spans the mathematical
frame of this paper. A \textit{tensor norm} $\alpha $ is a mapping which
assigns to each pair $(X,Y)$ of Banach spaces a norm $\alpha (\cdot ;X,Y)$
on the algebraic tensor product $X\otimes Y$ (shorthand: $X$ ${\otimes }_{%
\stackrel{}{\alpha }}$ $Y$ and $X\tilde{\otimes}_{\stackrel{}{\alpha }}Y$
for the completion) so that

\begin{itemize}
\item  $\varepsilon \leq \alpha \leq \pi $

\item  $\alpha $ satisfies the metric mapping property: If $S\in {\frak{L}}%
(X,Z)$ and $T\in {\frak{L}}(Y,W)$, then $\Vert S\otimes T:X\otimes _{%
\stackrel{}{\alpha }}Y\longrightarrow Z\otimes _{\stackrel{}{\alpha }}W\Vert
\leq \Vert S\Vert \ \Vert T\Vert $ .
\end{itemize}

\noindent Wellknown examples are the injective tensor norm $\varepsilon $,
which is the smallest one, and the projective tensor norm $\pi $, which is
the largest one. For other important examples we refer to \cite{df}, \cite
{gl}, or \cite{l}. Each tensor norm $\alpha $ can be extended in two natural
ways. For this, denote for given Banach spaces $X$ and $Y$
\[
\text{FIN}(X):=\{E\subseteq X\mid E\in \text{FIN}\}\hspace{0.2cm}\text{and}%
\hspace{0.2cm}\text{COFIN}(X):=\{L\subseteq X\mid X/L\in \text{FIN}\}\text{,}
\]
where FIN stands for the class of all finite dimensional Banach spaces. Let $%
z\in X\otimes Y$. Then the \textit{finite hull}\/ $\stackrel{\rightarrow }{%
\alpha }$ is given by
\[
\stackrel{\rightarrow }{\alpha }(z;X,Y):=\inf \{\alpha (z;E,F)\mid E\in
\text{FIN}(X)\text{, }F\in \text{FIN}(Y)\text{, }z\in E\otimes F\}\text{,}
\]
and the \textit{cofinite hull}\/ $\stackrel{\leftarrow }{\alpha }$ of $%
\alpha $ is given by
\[
\hspace{0.2cm}\stackrel{\leftarrow }{\alpha }(z;X,Y):=\sup \{\alpha
(Q_K^X\otimes Q_L^Y(z);X/K,Y/L)\mid K\in \text{COFIN}(X)\text{, }L\in \text{
COFIN}(Y)\}\text{.}
\]
$\alpha $ is called \textit{finitely generated}\/ if $\alpha $ $=$ $%
\stackrel{\rightarrow }{\alpha }$, \textit{cofinitely generated}\/ if $%
\alpha $ $=$ $\stackrel{\leftarrow }{\alpha }$ (it is always true that $%
\stackrel{\leftarrow }{\alpha }$ $\leq \alpha $ $\leq $ $\stackrel{%
\rightarrow }{\alpha }$). $\alpha $ is called \textit{right--accessible} if $%
\stackrel{\leftarrow }{\alpha }$$(z;E,Y)$ $=$ $\stackrel{\rightarrow }{%
\alpha }(z;E,Y)$ for all $(E,Y)\in $ FIN $\times $ BAN, \textit{%
left--accessible} if $\stackrel{\leftarrow }{\alpha }$$(z;X,F)$ $=$ $%
\stackrel{\rightarrow }{\alpha }$$(z;X,F)$ for all $(X,F)\in $ BAN $\times $
FIN, and \textit{accessible} if it is right--accessible and
left--accessible. $\alpha $ is called \textit{totally accessible} if $%
\stackrel{\leftarrow }{\alpha }$ $=$ $\stackrel{\rightarrow }{\alpha }$. The
injective norm $\varepsilon $ is totally accessible, the projective norm $%
\pi $ is accessible -- but not totally accessible, and Pisier's construction
implies the existence of a finitely generated tensor norm which is neither
left-- nor right--accessible (see \cite[31.6.]{df}).

There exists a powerful one--to--one correspondence between finitely
generated tensor norms and maximal Banach ideals which links thinking in
terms of operators with ''tensorial'' thinking and which allows to transfer
notions in the ''tensor language'' to the ''operator language'' and
conversely. In particular, this one--to--one correspondence helps to extend
the trace duality
\[
<S,T>:=tr(TS)\;\;(S\in {\frak{F}}(X,Y),\;T\in {\frak{F}}(Y,X))
\]
to operator ideals by using tensor product methods. We refer the reader to
\cite{df} and \cite{oe1} for detailed informations concerning this subject.
Let $X,Y$ be Banach spaces and $z=\sum\limits_{i=1}^nx_i^{\prime }\otimes
y_i $\ be an Element in $X^{\prime }\otimes Y$. Then $T_z(x):=\sum%
\limits_{i=1}^n\langle x,x_i^{\prime }\rangle $ $y_i$ defines a
finite rank operator $T_z\in {\frak{F}}(X,Y)$ which is
independent of the representation of $z$ in $X^{\prime }\otimes
Y$. Let $\alpha $ be a finitely generated tensor norm and
$(\frak{A},\mathbf{A})$ be a maximal Banach ideal. $\alpha $ and
$(\frak{A},\mathbf{A})$ are said to be \textit{associated},
notation:
\[
(\frak{A},\mathbf{A})\sim \alpha \hspace{0.15cm}\text{(shorthand: }\frak{A}%
\sim \alpha \text{, resp.}\hspace{0.1cm}\alpha \sim \frak{A}\text{),}
\]
if for all $E,F\in $ FIN
\[
{\frak{A}}(E,F)=E^{\prime }{\otimes }_{\stackrel{}{\alpha }}F
\]
holds isometrically: ${\mathbf{A}}(T_z)=\alpha (z;E^{\prime },F)$.

Since we will use them throughout in this paper, let us recall the important
notions of the conjugate operator ideal (cf. \cite{glr}, \cite{jo} and \cite
{oe2}) and the adjoint operator ideal (all details can be found in the
standard references \cite{df} and \cite{p1}). Let $(\frak{A},\mathbf{A})$ be
a quasi--Banach ideal.

\begin{itemize}
\item  Let ${\frak{A}}^\Delta (X,Y)$ be the set of all $T\in {\frak{L}}(X,Y)$
which satisfy
\[
{\mathbf{A}}^\Delta (T):=\sup \{\mid tr(TL)\mid \text{ }\mid L\in {\frak{F}}%
(Y,X),{\mathbf{A}}(L)\leq 1\}<\infty .
\]
Then a Banach ideal $({\frak{A}}^\Delta ,{\mathbf{A}}^\Delta )$ is
obtained (here, $tr(\cdot )$ denotes the usual trace for finite
rank operators). It is called the \textit{conjugate ideal} of
$(\frak{A},\mathbf{A})$.

\item  Let ${\frak{A}}^{*}(X,Y)$ be the set of all $T\in {\frak{L}}(X,Y)$
which satisfy
\[
{\mathbf{A}}^{*}(T):=\sup \{\mid tr(TJ_E^XSQ_K^Y\}\mid \text{
}\mid E\in \text{FIN}(X),K\in \text{COFIN}(Y),{\mathbf{A}}(S)\leq 1\}<\infty .
\]
Then a Banach ideal $({\frak{A}}^{*},{\mathbf{A}}^{*})$ is
obtained. It is called the \textit{adjoint operator ideal }of
$(\frak{A},\mathbf{A})$.
\end{itemize}

By definition, it immediately follows that ${\frak{A}}^\Delta \stackrel{1}{%
\subseteq }{\frak{A}}^{*}$. Another easy, yet important
observation is the
following: let $(\frak{A},\mathbf{A})$ be a quasi--Banach ideal and $(\frak{B%
},\mathbf{B})$ be a quasi--Banach ideal. If $\frak{A}\stackrel{\frak{E}}{%
\subseteq }\frak{B}$, then ${\frak{B}}^{*}\stackrel{1}{\subseteq }{\frak{A}}^{*}$%
, and $\frak{A}\stackrel{\frak{F}}{\subseteq }\frak{B}$ implies
the inclusion ${\frak{B}}^\Delta \stackrel{1}{\subseteq
}{\frak{A}}^\Delta $. In particular, it follows that
${\frak{A}}^{\Delta *}\stackrel{1}{=}{\frak{A}}^{**}$ and
$({\frak{A}}^{\Delta \Delta })^{*}\stackrel{1}{=}{\frak{A}}^{*}$.
A deeper
investigation of relations between the Banach ideals $({\frak{A}}^\Delta ,%
{\mathbf{A}}^\Delta )$ and $({\frak{A}}^{*},{\mathbf{A}}^{*})$
needs the help of an important local property, known as
accessibility, which can be viewed as a local version of
injectivity and surjectivity. All necesary details about
accessibility of operator ideals and its applications can be
found in \cite {df}, \cite{oe2}, \cite{oe3} and \cite{oe4}. So
let us recall :

\begin{itemize}
\item  A quasi--Banach ideal $(\frak{A},\mathbf{A})$ is called \textit{%
right--accessible}, if for all $(E,Y)\in $ FIN $\times $ BAN, operators $%
T\in {\frak{L}}(E,Y)$ and $\varepsilon >0$ there are $F\in $ FIN$(Y)$ and $%
S\in {\frak{L}}(E,F)$ so that $T=J_F^YS$ and ${\mathbf{A}}(S)\leq
(1+\varepsilon ){\mathbf{A}}(T)$.

\item  $(\frak{A},\mathbf{A})$ is called \textit{left--accessible}, if for
all $(X,F)\in $ BAN $\times $ FIN, operators $T\in {\frak{L}}(X,F)$ and $%
\varepsilon >0$ there are $L\in $ COFIN$(X)$ and $S\in
{\frak{L}}(X/L,F)$ so that $T=SQ_L^X$ and ${\mathbf{A}}(S)\leq
(1+\varepsilon ){\mathbf{A}}(T)$.

\item  A left--accessible and right--accessible quasi--Banach ideal is called%
{\ }\textit{accessible}.

\item  $(\frak{A},\mathbf{A})$ is \textit{totally accessible}, if for every
finite rank operator $T\in {\frak{F}}(X,Y)$ acting between Banach spaces $X$%
, $Y$ and $\varepsilon >0$ there are $(L,F)\in $ COFIN$(X)\times
$ FIN$(Y)$ and $S\in {\frak{L}}(X/L,F)$ so that $T=J_F^YSQ_L^X$
and ${\mathbf{A}}(S)\leq (1+\varepsilon ){\mathbf{A}}(T)$.
\end{itemize}

Given quasi--Banach ideals $(\frak{A},\mathbf{A)}$ and
$(\frak{B},\mathbf{B)} $, let $({\frak{A}}\circ
{\frak{B}},\mathbf{A\circ B})$ be the corresponding product ideal
and $({\frak{A}}\circ {\frak{B}}^{-1},{\mathbf{A}}\circ {\mathbf{B}}^{-1})$
(resp. $({\frak{A}}^{-1}\circ
{\frak{B}},{\mathbf{A}}^{-1}\circ{\mathbf{B}}))$ the
corresponding ''right--quotient'' (resp. ''left--quotient''). We write $(%
{\frak{A}}^{inj},{\mathbf{A}}^{inj})$, to denote the \textit{injective hull}%
\emph{\ }of $\frak{A}$, the unique smallest injective quasi--Banach ideal
which contains $({\frak{A}},\mathbf{A})$, and $({\frak{A}}^{sur},{\mathbf{A}}%
^{sur})$, the \textit{surjective hull} of $\frak{A}$, is the
unique smallest surjective quasi--Banach ideal which contains
$({\frak{A}},\mathbf{A})$. Of particular importance are the
quotients ${\frak{A}}^{\dashv }:={\frak{I}}\circ {\frak{A}}^{-1}$
and ${\frak{A}}^{\vdash }:={\frak{A}}^{-1}\circ {\frak{I}}$ and
their relations to ${\frak{A}}^\Delta $ and ${\frak{A}}^{*}$,
treated in detail in \cite{oe1} and \cite{oe4}.

In addition to the maximal Banach ideal
$({\frak{L}},\mathcal{\Vert \cdot \Vert })\sim
\mathcal{\varepsilon }$ we mainly will be concerned with the
maximal Banach ideals $({\frak{I}},\mathbf{I})\sim
\mathcal{\mathbf{\pi }}$ (integral operators),
$({\frak{L}}_2,{\mathbf{L}}_2)\sim w_2$ (Hilbertian
operators), $({\frak{D}}_2,{\mathbf{D}}_2)\stackrel{1}{=}
({\frak{L}}_2^{*},{\mathbf{%
L}}_2^{*})\stackrel{1}{=}{\frak{P}}_2^d\circ {\frak{P}}_2
\sim w_2^{*}$\/ ($2$--dominated operators), $({\frak{P}}_p,{\mathbf{P}}_p)
\sim g_p\backslash =g_q^{*} $ (absolutely $p$--summing operators),
$1\leq p\leq \infty ,\frac
1p+\frac 1q=1$, $({\frak{L}}_\infty ,{\mathbf{L}}_\infty )\stackrel{1}{=}
({\frak{A}}_1^{*},{\mathbf{P}}_1^{*})\sim w_\infty $ and
$({\frak{L}}_1,{\mathbf{L}}_1)%
\stackrel{1}{=}({\frak{P}}_1^{*d},{\mathbf{P}}_1^{*d})\sim w_1$. We
also consider the maximal Banach ideals
$({\frak{C}}_2,{\mathbf{C}}_2)$ $\sim c_2$ (cotype 2 operators) and
$({\frak{A}}_P,{\mathbf{A}}_P)$ $\sim \alpha _P$ (Pisier`s
counterexample of a maximal Banach ideal which is neither right--
nor left--accessible (cf. \cite{df}, 31.6)).

\section{The principle of local reflexivity for operator ideals}

Let $(\frak{A},\mathbf{A})$ be a \textit{maximal} Banach ideal.
Then, ${\frak{A}}^\Delta $ always is right--accessible (cf.
\cite{oe4}). The natural question whether ${\frak{A}}^\Delta $ is
\textit{left}--accessible is still
open$\footnote{%
For minimal Banach ideals $(\frak{A},\mathbf{A)}$, there exist
counterexamples: the conjugate of ${\frak{A}}_P^{\min }$ neither
is right--accessible nor left--accessible (cf. \cite{oe2}).}$ and
leads to
interesting and non--trivial results concerning the local structure of $%
{\frak{A}}^\Delta $. Deeper investigations of the left--accessibility of $%
{\frak{A}}^\Delta $ namely lead to a link with a principle of
local reflexivity for operator ideals (a detailed discussion can
be found in \cite {oe1} and \cite{oe2}) which allows a
transmission of the operator norm estimation in the classical
principle of local reflexivity to the ideal norm $\mathbf{A}$. So
let us recall the

\begin{definition}
Let $E$ and $Y$ be Banach spaces, $E$ finite dimensional, $F\in $ FIN$%
(Y^{\prime })$ and $T\in {\frak{L}}(E,Y^{\prime \prime })$. Let $(\frak{A},%
\mathbf{A})$ be a quasi--Banach ideal and $\epsilon >0$. We say
that the principle of ${\frak{A}}-$local reflexivity (short:
${\frak{A}}-LRP$) is satisfied, if there exists an operator $S\in
{\frak{L}}(E,Y)$ so that

\begin{enumerate}
\item[(1)]  ${\mathbf{A}}(S)\leq (1+\epsilon )\cdot {\mathbf{A}}^{**}(T)$

\item[(2)]  $\left\langle Sx,y^{\prime }\right\rangle =\left\langle
y^{\prime },Tx\right\rangle $ for all $(x,y^{\prime })\in E\times F$

\item[(3)]  $j_YSx=Tx$ for all $x\in T^{-1}(j_Y(Y))$.
\end{enumerate}
\end{definition}

Although both, the quasi--Banach ideal $\frak{A}$ and\textit{\ }the $1$%
--Banach ideal ${\frak{A}}^{**}$ are involved, the unbalance can
be justified by the following statement which holds for arbitrary
quasi--Banach ideals (see \cite{oe2}):

\begin{theorem}
Let $(\frak{A},\mathbf{A})$ be a quasi--Banach ideal. Then the following
statements are equivalent:

\begin{enumerate}
\item[(i)]  ${\frak{A}}^\Delta $ is left--accessible

\item[(ii)]  ${\frak{A}}^{**}(E,Y^{\prime \prime })$ $\widetilde{=}$ ${\frak{A}}%
(E,Y)^{\prime \prime }$ for all $(E,Y)\in $ FIN $\times $ BAN

\item[(iii)]  The ${\frak{A}}-LRP$ holds.
\end{enumerate}
\end{theorem}

If we only assume that $(\frak{A},\mathbf{A})$ is a maximal Banach ideal,
the previous statement is contained in the following important observation
(cf. \cite{oe1}):

\begin{theorem}
Let $(\frak{A},\mathbf{A})\sim \alpha $ be an arbitrary maximal Banach
ideal. Then the following statements are equivalent:

\begin{enumerate}
\item[(i)]  ${\frak{A}}^\Delta $ is left--accessible

\item[(ii)]  ${\frak{A}}(E,Y^{\prime \prime })$ $\widetilde{=}$ ${\frak{A}}%
(E,Y)^{\prime \prime }$ for all $(E,Y)\in $ FIN $\times $ BAN

\item[(iii)]  $({\frak{A}}^d)^\Delta (X_0,Y^{\prime })\;\widetilde{=}%
\;(X_0\otimes _{\overleftarrow{\alpha }}Y)^{\prime }$ for all Banach spaces $%
X_0$ with the metric approximation property and all $Y\in $ BAN

\item[(iv)]  The ${\frak{A}}-LRP$ holds.
\end{enumerate}
\end{theorem}

One reason which leads to extreme persistent difficulties
concerning the verification of the ${\frak{A}}-LRP$ for a given
maximal Banach ideal $\frak{A} $, is the behaviour of the bidual
$({\frak{A}}^\Delta )^{dd}$: although we know that in general
$({\frak{A}}^\Delta )^{dd}$ and $({\frak{A}}^\Delta )^d$
both are accessible (see \cite{oe1} and \cite{oe2}) and that $({\frak{A}}%
^\Delta )^{dd}\stackrel{1}{\subseteq }{\frak{A}}^\Delta $, we do
not know whether ${\frak{A}}^\Delta (X,Y)$ and $({\frak{A}}^\Delta
)^{dd}(X,Y)$ coincide \textit{isometrically} for \textit{all}
Banach spaces $X$ and $Y$. If we allow in addition the
approximation property of $X$ or $Y$, then we may state the
following

\begin{theorem}
Let $(\frak{A},\mathbf{A})$ be an arbitrary maximal Banach ideal and $X$, $Y$
be arbitrary Banach spaces. Then
\[
{\frak{A}}^{d\Delta }(X,Y)\stackrel{1}{=}{\frak{A}}^{\Delta
d}(X,Y)
\]
holds in each of the following two cases:

\begin{enumerate}
\item[(i)]  $X^{\prime }$ has the metric approximation property

\item[(ii)]  $Y^{\prime }$ has the metric approximation property and the $%
{\frak{A}}^d-LRP$ is satisfied.
\end{enumerate}
\end{theorem}

Which class of operator ideals $(\frak{A},\mathbf{A})$ does satisfy the $%
{\frak{A}}-LRP$? At least, all maximal and right--accessible
Banach ideals belong to this class (cf. \cite{oe5}) due to the
following

\begin{theorem}
Let $(\frak{A},\mathbf{A)}$ be an arbitrary Banach ideal. If
$\frak{A}$ is right--accessible and ultrastable, then the
${\frak{A}}-LRP$ is satisfied.
\end{theorem}

Pisier`s counterexample of the maximal Banach ideal $({\frak{A}}_P,{\mathbf{A}}%
_P)$ which neither is left--accessible nor right--accessible (cf.
\cite{df}, 31.6) implies that in particular
$({\frak{A}}_P^{*},{\mathbf{A}}_P^{*})$ neither is
left--accessible nor right--accessible. Thinking at
${\frak{A}}_P^\Delta \stackrel{1}{\subseteq }{\frak{A}}_P^{*}$,
this leads to the natural and even more tough question whether the
${\frak{A}}_P-LRP$ is true or false. Since the
Pisier space $P$ does not have the approximation property, ${\frak{A}}%
_P^\Delta $ cannot be totally\textit{\ }accessible. Is it even true that $(%
{\frak{A}}_P^\Delta )^{inj}\stackrel{1}{=}{\frak{P}}_1\circ
({\frak{A}}_P)^{-1}$ is not totally accessible? If this is the
case, the ${\frak{A}}_P-LRP$ will be false. However,
${\frak{A}}_P^\Delta $ is not injective (cf. \cite{oe5}). What
about the left accessibility of ${\frak{A}}_P^{*\Delta }$? Can we
mimic Pisier`s proof to construct a similar counterexample of a
maximal Banach
ideal $({\frak{A}}_0,{\mathbf{A}}_0)$ which even does not satisfy the ${\frak{A}}%
_0-LRP\;$(cf. \cite[31.6.]{df})? Unfortunately, the proof of the following
statement leads to a negative answer:

\begin{proposition}
Let $Y$ be an arbitrary Banach space. Then
\[
({\frak{A}}_P^{*\Delta })^{inj}(P,Y)\subseteq
{\frak{L}}_2(P,Y)\text{,}
\]
and
\[
{\mathbf{L}}_2(T)\leq (2{\mathbf{C}}_2(P^{\prime })\cdot
{\mathbf{C}}_2(P))^{\frac 32}\cdot ({\mathbf{A}}_P^{*\Delta
})^{inj}(T)
\]
for all operators $T\in ({\frak{A}}_P^{*\Delta })^{inj}(P,Y)$.
\end{proposition}

\textsc{Proof:} First, let $F$ $\in $ FIN and $T\in
{\frak{L}}(P,F)$ an
arbitrary linear operator. Since the bidual $(({\frak{A}}_P^{*\Delta })^{dd},(%
{\mathbf{A}}_P^{*\Delta })^{dd})=:({\frak{B}}_P,{\mathbf{B}}_P)$
always is left--accessible (!) and $T$ is a finite rank operator,
a copy of Pisier`s construction immediately leads to
\[
{\mathbf{L}}_2(T)\leq (2{\mathbf{C}}_2(P^{\prime })\cdot
{\mathbf{C}}_2(P))^{\frac
32}\cdot {\mathbf{A}}_P^{*\Delta }(T^{\prime \prime })=(2{\mathbf{C}}%
_2(P^{\prime })\cdot {\mathbf{C}}_2(P))^{\frac 32}\cdot
{\mathbf{B}}_P(T)\text{.}
\]
Conjugation therefore implies
\[
{\frak{D}}_2(F,P)\stackrel{1}{=}{\frak{L}}_2^\Delta (F,P)\subseteq {\frak{B}}%
_P^\Delta (F,P)\stackrel{1}{\subseteq }{\frak{A}}_P^{*}(F,P)
\]
and
\[
{\mathbf{A}}_P^{*}(S)\leq (2{\mathbf{C}}_2(P^{\prime })\cdot {\mathbf{C}}%
_2(P))^{\frac 32}\cdot {\mathbf{D}}_2(S)\text{ \ for all }S\in {\frak{L}}(F,P)%
\text{.}
\]
Since $({\frak{A}}_P^{*},{\mathbf{A}}_P^{*})$ is maximal, it even
follows that
\begin{equation}
{\frak{D}}_2(Y,P)\subseteq {\frak{A}}_P^{*}(Y,P)
%\tag{$*$}
\end{equation}
for all Banach spaces $Y$ and
\[
{\mathbf{A}}_P^{*}(S)\leq (2{\mathbf{C}}_2(P^{\prime })\cdot {\mathbf{C}}%
_2(P))^{\frac 32}\cdot {\mathbf{D}}_2(S)\text{ \ for all \ }S\in {\frak{L}}(Y,P)%
\text{.}
\]
Hence, conjugation of $(1)$ finishes the proof.$\blacksquare $

Note that the accessibility of the bidual $({\frak{A}}^\Delta
)^{dd}$ (which also was used in the previous proof) implies one
of the main difficulties which appear repeatedly if one tries to
construct a counterexample of a maximal Banach ideal
$(\frak{A},\mathbf{A})$ so that the ${\frak{A}}-LRP$ is not
satisfied. In general, one is allowed to substitute statements
related to properties of ${\frak{A}}^\Delta $ through statements
related to properties of the (left--)accessible bidual
$({\frak{A}}^\Delta )^{dd}$ so that these statements remain to be
true, regardless whether the ${\frak{A}}-LRP$ is satisfied or
not! In particular, such statements cannot be used for a proof by
contradiction. However, a first step towards a successful
construction of such a candidate $(\frak{A},\mathbf{A})$ is given
by the following factorization property for finite rank operators
which had been introduced by Jarchow and Ott (cf. \cite{jo}).
This factorization property not only turns out to be a useful
tool in constructing such a counterexample; it even allows one to
show that ${\frak{L}}_\infty $ \textit{is not totally accessible}
-- solving a problem of Defant and Floret (see \cite[21.12.]{df}
and \cite[Theorem 4.1]{oe5}). So, let us recall the definition of
this factorization property and its implications:

\begin{definition}[Jarchow/Ott]
Let $(\frak{A},\mathbf{A)}$ and $(\frak{B},\mathbf{B)}$ be
arbitrary quasi--Banach ideals. Let $L\in {\frak{F}}(X,Y)$ an
arbitrary finite rank operator between two Banach spaces $X$ and
$Y$. Given $\epsilon >0$, we can
find a Banach space $Z$ and operators $A\in {\frak{A}}(Z,Y)$, $B\in {\frak{B}}%
(X,Z)$ so that $L=AB$ and
\[
{\mathbf{A}}(A)\cdot {\mathbf{B}}(B)\leq (1+\epsilon )\cdot {\mathbf{A\circ B}}(L)%
\text{.}
\]

\begin{enumerate}
\item[(i)]  If the operator $A$ is of finite rank, we say that $\frak{A}%
\circ \frak{B}$ has the property\textit{\ }(I).

\item[(ii)]  If the operator $B$ is of finite rank, we say that $\frak{A}%
\circ \frak{B}$ has the property\textit{\ }(S).
\end{enumerate}
\end{definition}

Important examples are the following (see \cite[Lemma 2.4.]{jo}):

\begin{itemize}
\item  If $\frak{B}$ is injective, or if $\frak{A}$ contains ${\frak{L}}_2$ as
a factor, then $\frak{A}\circ \frak{B}$ has the property\textit{\ }(I).

\item  If $\frak{A}$ is surjective, or if $\frak{B}$ contains ${\frak{L}}_2$
as a factor, then $\frak{A}\circ \frak{B}$ has the property\textit{\ }(S).
\end{itemize}

Since ${\frak{L}}_2\circ \frak{A}$ is injective for every quasi--Banach ideal $%
(\frak{A},\mathbf{A})$ (see \cite[Lemma 5.1.]{oe4}), $\frak{B}\circ {\frak{L}}%
_2\circ \frak{A}$ therefore has the property (I) as well as the
property (S), for all quasi--Banach ideals
$(\frak{A},\mathbf{A})$ and $(\frak{B},\mathbf{B})$. Such ideals
are exactly those which contain ${\frak{L}}_2$ as factor -- in
the sense of \cite {jo}. \cite{oe5} explains in detail how the
property (I) influences the
structure of operator ideals of type ${\frak{A}}^{inj*}\stackrel{1}{=}%
\diagdown {\frak{A}}^{*}$ and their conjugates, leading to
another approach to construct a counterexample of a maximal
Banach ideal with non--left--accessible conjugate. To this end,
first note that for all Banach spaces $X$,$Y$ and
$X\stackrel{1}{\hookrightarrow }Z$, every operator $T\in
({\frak{A}}^{inj})^{*}(X,Y)\stackrel{1}{=}\diagdown
{\frak{A}}^{*}(X,Y)$ satisfies the following extension property:
given $\epsilon >0$, there exists an operator $\widetilde{T}\in
\diagdown {\frak{A}}^{*}(Z,Y^{\prime
\prime })$ so that $j_YT=\widetilde{T}J_X^Z$ and $\diagdown {\mathbf{A}}^{*}(%
\widetilde{T})\leq (1+\epsilon )\cdot \diagdown
{\mathbf{A}}^{*}(T)$ (see \cite[Satz 7.14]{h1}). In particular,
such an extension holds for all finite rank operators. However,
we then cannot be sure that $\widetilde{T}$ is also as a\textit{\
finite rank }operator. Here, property (I) comes into play -- in
the following sense:

\begin{theorem}
Let $(\frak{A},\mathbf{A)}$ be a Banach ideal so that
${\frak{A}}^{*}\circ {\frak{L}}_\infty $ has the property\textit{\
}(I). Let $\epsilon >0$, $X$ and $Y$ be arbitrary Banach spaces
and $L\in {\frak{F}}(Y,X)$. Let $Z$ be a Banach space which
contains $Y$ isometrically. Then there exists a finite rank
operator $V\in {\frak{F}}(Z,X^{\prime \prime })$ so that
$j_XL=VJ_Y^Z$ and
\[
({\mathbf{A}}^{inj})^{*}(V)\leq (1+\epsilon )\cdot
({\mathbf{A}}^{inj})^{*}(L)\text{.}
\]
If in addition the ${\frak{A}}^{*}-LRP$ is satisfied, then $V$
even can be chosen to be a finite rank operator with range in $X$
and $L=VJ_Y^Z$.
\end{theorem}

Consequently, this theorem leads to important implications which link the
principle of local reflexivity for operator ideals with product ideals of
type (I) such as the following ones (cf. \cite{oe5}):

\begin{theorem}
Let $(\frak{A},\mathbf{A})$ be a Banach ideal so that the
${\frak{A}}^{*}-LRP$ is satisfied. Then
\[
{\frak{A}}^{*\Delta inj}\stackrel{1}{\subseteq
}({\frak{A}}^{*\Delta inj})^{dd}
\]
If in addition, ${\frak{A}}^{*}\circ {\frak{L}}_\infty $ has the
property (I), then
\[
({\frak{A}}^{*\Delta inj})^{dd}\stackrel{1}{=}{\frak{A}}^{*\Delta inj}\stackrel{1%
}{=}{\frak{P}}_1\circ ({\frak{A}}^{*})^{-1}\stackrel{1}{=}{\frak{A}}^{inj*\Delta }%
\stackrel{1}{=}({\frak{A}}^{inj*\Delta })^{dd}
\]
\end{theorem}

\begin{theorem}
Let $(\frak{A},\mathbf{A})$ be a left--accessible Banach ideal so that $%
{\frak{A}}^{*}\circ {\frak{L}}_\infty $ has the property (I). Then ${\frak{A}}%
^{inj}$ is totally accessible and ${\frak{A}}^{inj}\stackrel{1}{\subseteq }(%
{\frak{A}}^{inj})^{*\Delta }$. If in addition
$(\frak{A},\mathbf{A})$ is maximal, then $({\frak{A}}^{inj})^{*}$
is also totally accessible.
\end{theorem}

\begin{theorem}
Let $(\frak{A},\mathbf{A})$ be a Banach ideal so that
${\frak{A}}^{*}\circ {\frak{L}}_\infty $ has the property (I). If
space$(\frak{A})$ contains a
Banach space $X_0$ so that $X_0$ has the bounded approximation property but $%
X_0^{\prime \prime }$ has not, then the ${\frak{A}}^{*}-LRP$
cannot be satisfied.
\end{theorem}

To construct such maximal Banach ideals, note again that
${\frak{A}}^{*}\circ
{\frak{L}}_\infty $ has the property (I), if ${\frak{A}}^{*}$ contains ${\frak{L}}%
_2$ as a factor. Since ${\frak{A}}^{*}$ is a Banach ideal, we
therefore have
to look for maximal \textit{Banach} ideals of type $\frak{B}\circ {\frak{L}}%
_2\circ \frak{C}$. A first investigation of geometrical properties of such
product ideals was given in \cite{oe4}.

\section{Normed operator ideal products}

Unfortunately, we still cannot present explicite sufficient criteria which
show the existence of (an equivalent) ideal norm on product ideals. It seems
to be much more easier to show that a certain product ideal cannot be a
normed one by using arguments which involve trace ideals and the ideal of
nuclear operators (the smallest Banach ideal). Even more holds: if $\frak{A}%
\circ {\frak{L}}_2$ is a 1--Banach ideal for certain operator ideals $\frak{A}$%
, then $\frak{A}\circ {\frak{L}}_2$ \textit{is not}
right--accessible (cf. Theorem 4.4)! However, let us look more
carefully at such product ideals. First, we note an improvement
of our own work (cf. \cite[Theorem 4.2]{oe5}):

\begin{theorem}
Let $(\frak{A},\mathbf{A})$ be a maximal Banach ideal. Then both,
the maximal $\frac 12$--Banach ideal ${\frak{A}}^{inj}\circ
{\frak{L}}_2$ and the injective hull of the maximal $\frac
12$--Banach ideal $\frak{A}\circ {\frak{A}}_2$ are totally
accessible.
\end{theorem}

\textsc{Proof:} Since every Hilbert space has the metric
approximation property and since
${\frak{A}}^{inj}\stackrel{1}{=}({\frak{A}}^{inj})^{**}$ is
right--accessible, an easy calculation shows that
\begin{equation}
{\frak{A}}^{inj}\circ
{\frak{L}}_2\stackrel{1}{=}({\frak{A}}^{inj})^{*\Delta }\circ
{\frak{L}}_2\text{.}
%\tag{$*$}
\end{equation}
Since $({\frak{A}}^{inj})^{*\Delta }$ is right--accessible, the
total
accessibility of ${\frak{L}}_2$ and the property (S) of the product ideal $(%
{\frak{A}}^{inj})^{*\Delta }\circ {\frak{L}}_2$ even imply that $({\frak{A}}%
^{inj})^{*\Delta }\circ {\frak{L}}_2$ is totally accessible (cf.
\cite[Proposition 4.1]{oe5}). Hence, ${\frak{A}}^{inj}\circ
{\frak{L}}_2$ is totally accessible (due to $(2)$), and in
particular we obtain that $({\frak{A}}\circ
{\frak{L}}_2)^{inj}\stackrel{1}{=}({\frak{A}}^{inj}\circ
{\frak{L}}_2)^{inj} $ is totally accessible.$\blacksquare $

Now, let $(\frak{A},\mathbf{A})$ be a maximal Banach ideal so that
${\mathbf{L}}_2\circ \mathbf{A}$ even is a norm on the (maximal) product ideal $({\frak{L}}%
_2\circ {\frak{A}},{\mathbf{L}}_2\circ {\mathbf{A}})$. Then ${\frak{A}}^{*}\stackrel{%
1}{\subseteq }({\frak{L}}_2\circ {\frak{A}})^{*}\stackrel{1}{\subseteq }{\frak{L}}%
_\infty $ (cf. \cite[Proposition 5.1.]{oe4}) and ${\frak{L}}_\infty \stackrel{1%
}{=}{\frak{P}}_1^\Delta \stackrel{1}{\subseteq }{\frak{N}}^\Delta
$. Given Banach spaces $X$ and $Y$ so that both, $X^{\prime }$
and $Y$ have cotype 2, \cite[Theorem 4.9.]{pi} tells us, that any
finite rank operator $L\in {\frak{F}}(Y,X)$ satisfies
\[
{\mathbf{N}}(L)\leq (2{\mathbf{C}}_2(X^{\prime })\cdot
{\mathbf{C}}_2(Y))^{\frac 32}\cdot {\mathbf{D}}_2(L)\text{.}
\]
Hence,
\[
{\frak{N}}^\Delta (X,Y)\subseteq {\frak{D}}_2^\Delta (X,Y)\stackrel{1}{=}{\frak{L}}%
_2(X,Y)\text{,}
\]
and we have proven a rather surprising fact (revealing the strong influence
of a \textit{norm} on an operator ideal product):

\begin{theorem}
Let $(\frak{A},\mathbf{A})$ be a maximal Banach ideal so that the
product ideal $({\frak{L}}_2\circ {\frak{A}},{\mathbf{L}}_2\circ
\mathbf{A})$ is normed. Let $X$ and $Y$ be arbitrary Banach
spaces so that both, $X^{\prime }$ and $Y $ have cotype 2. Then
\[
{\frak{A}}^{*}(X,Y)\stackrel{1}{\subseteq }({\frak{L}}_2\circ {\frak{A}}%
)^{*}(X,Y)\subseteq {\frak{L}}_2(X,Y)
\]
and
\[
{\mathbf{L}}_2(T)\leq (2{\mathbf{C}}_2(X^{\prime })\cdot
{\mathbf{C}}_2(Y))^{\frac 32}\cdot ({\mathbf{L}}_2\circ
{\mathbf{A}})^{*}(T)\leq (2{\mathbf{C}}_2(X^{\prime })\cdot
{\mathbf{C}}_2(Y))^{\frac 32}\cdot {\mathbf{A}}^{*}(T)
\]
for all operators $T\in {\frak{A}}^{*}(X,Y)$.
\end{theorem}

To maintain the previous statement, even a permutation of the factors $\frak{%
A}$ and ${\frak{L}}_2$ in the product ${{\frak{L}}}_2\circ
\frak{A}$ is allowed:

\begin{theorem}
Let $(\frak{A},\mathbf{A})$ be a maximal Banach ideal so that the
product ideal $({\frak{A}}\circ {\frak{L}}_2,{\mathbf{A\circ
L}}_2)$ is normed. Let $X$ and $Y$ be arbitrary Banach spaces so
that both, $X^{\prime }$ and $Y$ have cotype 2. Then
\[
{\frak{A}}^{*}(X,Y)\stackrel{1}{\subseteq }({\frak{A}}\circ {\frak{L}}%
_2)^{*}(X,Y)\subseteq {\frak{L}}_2(X,Y)
\]
and
\[
{\mathbf{L}}_2(T)\leq (2{\mathbf{C}}_2(X^{\prime })\cdot
{\mathbf{C}}_2(Y))^{\frac 32}\cdot ({\mathbf{A\circ
L}}_2{\mathbf{)}}^{*}(T)\leq (2{\mathbf{C}}_2(X^{\prime })\cdot
{\mathbf{C}}_2(Y))^{\frac 32}\cdot {\mathbf{A}}^{*}(T)
\]
for all operators $T\in {\frak{A}}^{*}(X,Y)$.
\end{theorem}

\textsc{Proof:} Let $(\frak{A},\mathbf{A})$ and $X$, $Y$ be as
before and let ${\frak{A}}\circ {\frak{L}}_2$ be normed. Then
${\frak{A}}\stackrel{1}{=}{\frak{A}}^{dd}$, and ${\frak{A}}\circ
{\frak{A}}_2$ is a maximal (and therefore a
regular) Banach ideal (cf. \cite[Lemma 4.3]{oe5}). Since the injective $%
\frac 12$--Banach ideal ${\frak{L}}_2\circ {\frak{A}}^d$ is also
regular (cf. \cite[Lemma 5.1]{oe4}), an easy calculation shows
that
\[
({\frak{A}}\circ {\frak{L}}_2)^d\stackrel{1}{=}{\frak{L}}_2\circ
{\frak{A}}^d
\]
and\footnote{%
In particular, it follows that ${\frak{A}}\circ {\frak{L}}_2$ is
surjective (cf. \cite[8.5.9.]{p1}).}
\[
{\frak{A}}\circ {\frak{L}}_2\stackrel{1}{=}({\frak{L}}_2\circ
{\frak{A}}^d)^d\text{.}
\]
Since ${\mathbf{A\circ L}}_2$ is a norm, $({\mathbf{A\circ L}}%
_2)^d$ obviously is a norm too. Hence, if we apply the previous theorem to
the normed product ideal $({\frak{A}}\circ {\frak{L}}_2)^d\stackrel{1}{=}{\frak{L}}%
_2\circ {\frak{A}}^d$, we obtain
\[
{\frak{A}}^{*d}(X,Y)\stackrel{1}{\subseteq }({\frak{L}}_2\circ {\frak{A}}%
^d)^{*}(X,Y)\stackrel{1}{=}({\frak{A}}\circ
{\frak{L}}_2)^{*d}(X,Y)\subseteq {\frak{L}}_2(X,Y)\text{,}
\]
and
\[
{\mathbf{L}}_2(T)\leq C\cdot ({\mathbf{A\circ
L}}_2)^{*}(T^{\prime })\leq C\cdot {\mathbf{A}}^{*}(T^{\prime })
\]
for \textit{all} operators $T\in {\frak{A}}^{*d}(X,Y)$ (where $C:=(2{\mathbf{C}}%
_2(X^{\prime })\cdot {\mathbf{C}}_2(Y))^{\frac 32}$). Now, since
$Y$ has the same coptype as its bidual $(Y^{\prime })^{\prime }$
with identical cotype
constants (cf. \cite[Corollary 11.9]{djt}), the proof is finished.$%
\blacksquare $

Let $(\frak{A},\mathbf{A})$ be a given ultrastable quasi--Banach
ideal so that $({\frak{A}}\circ {\frak{L}}_2,{\mathbf{A\circ
L}}_2)$ is right--accessible.
Our aim is to show that in this case $({\frak{A}}\circ {\frak{L}}_2,{\mathbf{%
A\circ L}}_2)$ and $({\frak{L}}_2\circ {\frak{A}}^{*},{\mathbf{L}}_2\circ {\mathbf{A}}%
^{*})$ \textit{both together} cannot be normed.\textit{\ }To this end, we
need a lemma which is of its own interest:

\begin{lemma}
Let $({\frak{A}}_0,{\mathbf{A}}_0)$ be a maximal Banach ideal so that space$(%
{\frak{A}}_0)$ contains a Banach space without the approximation
property. Then there does not exist a maximal Banach ideal
$(\frak{C},\mathbf{C})$ so
that ${\frak{C}}\circ {\frak{L}}_\infty $ has the property (I) and ${\frak{C}}%
\subseteq {\frak{A}}_0^{-1}\circ {\frak{P}}_1$.
\end{lemma}

\textsc{Proof:} Assume that the statement is false. Then there
exists a (maximal) Banach ideal $(\frak{A},\mathbf{A})$ so that
${{\frak{A}}}_0\subseteq {{\frak{P}}}_1\circ
({\frak{A}}^{*})^{-1}\stackrel{1}{=}({\frak{A}}^{*\Delta })^{inj}
$. Due to the assumed property (I) of ${\frak{A}}^{*}\circ
{\frak{L}}_\infty $, the proof of Theorem 3.4 in \cite{oe5} shows
that even $(({\frak{A}}^{*\Delta
})^{inj})^{dd}\stackrel{1}{\subseteq }({\frak{A}}^{inj})^{*\Delta }\stackrel{1%
}{\subseteq }{\frak{N}}^\Delta $. Since ${\frak{A}}_0$ was assumed
to be a
maximal Banach ideal, we therefore obtain ${\frak{A}}_0\stackrel{1}{=}{\frak{A}}%
_0^{dd}\stackrel{1}{\subseteq }{\frak{N}}^\Delta $ which is a contradiction.$%
\blacksquare $

\begin{corollary}
Let $({\frak{A}}_0,{\mathbf{A}}_0)$ be a maximal Banach ideal so that space$(%
{\frak{A}}_0)$ contains a Banach space without the approximation property. If $%
({\frak{A}}_0^{-1}\circ {\frak{P}}_1)\circ {\frak{L}}_\infty $ has
the property (I), ${\frak{A}}_0$ is not left--accessible.
\end{corollary}

\begin{theorem}
Let $(\frak{B},\mathbf{B})$ be an ultrastable quasi--Banach ideal so that $%
{\frak{B}}\subseteq {\frak{L}}_\infty $. If ${\frak{B}}\circ
{\frak{L}}_2$ is right--accessible, ${\frak{B}}\circ {\frak{L}}_2$
cannot be a $1-$ Banach ideal.
\end{theorem}

\textsc{Proof:} Assume that the statement is false and put ${\frak{B}}_0:=(%
{\frak{L}}_\infty \circ {\frak{L}}_2)^{*}$ and ${\frak{A}}:=({\frak{B}}\circ {\frak{L}}%
_2)^{*}$. Then
\[
{\frak{A}}^{*}\stackrel{1}{=}({\frak{B}}\circ {\frak{L}}_2)^{\max }\stackrel{1}{=}(%
{\frak{B}}\circ
{\frak{L}}_2)^{reg}\stackrel{1}{=}{\frak{B}}^{reg}\circ
{\frak{L}}_2
\]
is right--accessible (cf. \cite[Proposition 2.3]{oe5}) and contains ${\frak{L}}%
_2$ as a factor so that in particular ${\frak{A}}^{*}\circ
{\frak{L}}_\infty $ has the property (I). Since
${\frak{B}}\subseteq {\frak{L}}_\infty $,
\[
{\frak{B}}_0\circ {\frak{A}}^{*}\subseteq {\frak{A}}\circ {\frak{A}}^{*}\stackrel{1}{%
\subseteq }{\frak{I}}\stackrel{1}{\subseteq }{\frak{P}}_1\text{,}
\]
and it follows that ${\frak{A}}^{*}\subseteq
{\frak{B}}_0^{-1}\circ {\frak{P}}_1$. Since $Id_P\in {\frak{B}}_0$
(cf. \cite[Proposition 4.4]{oe5}), Lemma 4.1 leads to a
contradiction.$\blacksquare $

Now let us assume that $(\frak{B},\mathbf{B})$ is even is a \textit{maximal
Banach ideal} so that ${\frak{B}}\subseteq {\frak{L}}_\infty $. If ${\frak{B}}%
\circ {\frak{L}}_2$ were normed, then ${\frak{B}}\circ
{\frak{L}}_2$ would be a
maximal \textit{and }surjective Banach ideal, implying that ${\frak{P}}_1^d%
\stackrel{1}{=}{\frak{I}}^{sur}\stackrel{1}{=}({\frak{N}}^{max})^{sur}%
\stackrel{1}{\subseteq }{\frak{B}}\circ {\frak{L}}_2$. Hence,
\begin{equation}
({\frak{B}}\circ {\frak{L}}_2)^{*}\stackrel{1}{\subseteq }{\frak{P}}_1^{d*}%
\stackrel{1}{=}{\frak{L}}_1\text{.}
%\tag{$*$}
\end{equation}
Since ${\frak{B}}\subseteq {\frak{L}}_\infty $, it even follows that ${\frak{P}}%
_1\stackrel{1}{=}{\frak{L}}_\infty ^{*}\subseteq ({\frak{B}}\circ {\frak{L}}_2)^{*}%
\stackrel{1}{\subseteq }{\frak{L}}_1$ which is a contradiction
(cf. \cite[27.2.]{df}). So, in this case we obtain a stronger
result:

\begin{theorem}
Let $(\frak{B},\mathbf{B})$ be a maximal Banach ideal so that ${\frak{B}}%
\subseteq {\frak{L}}_\infty $. Then ${\frak{B}}\circ {\frak{L}}_2$
cannot be a $1-$ Banach ideal.
\end{theorem}

We finish this survey paper with two statements linking the principle of
local reflexivity with \textit{normed }operator ideal products consisting of
operators which factor through a Hilbert space (cf. \cite[Theorem 4.4]{oe5}):

\begin{theorem}
Let $(\frak{B},\mathbf{B})$ be an ultrastable quasi--Banach ideal and $X_0$
be a Banach space without the bounded approximation property so that
\[
({\frak{B}}\circ {\frak{L}}_2)^{reg}(\cdot ,X_0)\stackrel{\frak{F}}{\subseteq }%
{\frak{P}}_1(\cdot ,X_0)\text{ .}
\]
If ${\frak{B}}\circ {\frak{L}}_2$ is a $1$--Banach ideal, then the $({\frak{B}}%
\circ {\frak{L}}_2)^{**}-LRP$ cannot be satisfied.
\end{theorem}

\begin{theorem}
Let $(\frak{B},\mathbf{B})$ be a maximal Banach ideal so that ${\frak{P}}%
_1\circ {\frak{B}}^{-1}$ contains ${\frak{L}}_2$ as a factor and ${\frak{L}}%
_\infty \circ ({\frak{P}}_1\circ {\frak{B}}^{-1})$ is normed. If
there exists a
Banach space $X_0$ without the approximation property so that $X_0\in $ space%
$(\frak{B})$, then the ${\frak{B}}-LRP$ is false.
\end{theorem}

\textsc{Proof:} Assume that the statement is not true. Then ${\frak{A}}_0:=%
{\frak{P}}_1\circ {\frak{B}}^{-1}\stackrel{1}{=}({\frak{B}}^\Delta
)^{inj}$ is
totally accessible. Since ${\frak{A}}_0$ contains ${\frak{L}}_2$ as a factor, $%
{\frak{L}}_\infty \circ {\frak{A}}_0$ in particular has the
property (S), and it follows that ${\frak{L}}_\infty \circ
{\frak{A}}_0$ even is totally accessible (cf. \cite[Proposition
4.1]{oe5}). Hence,
\[
{\frak{B}}\stackrel{1}{\subseteq }{\frak{A}}_0^{-1}\circ {\frak{P}}_1\stackrel{1}{=%
}({\frak{L}}_\infty \circ
{\frak{A}}_0)^{*}\stackrel{1}{=}({\frak{L}}_\infty \circ
{\frak{A}}_0)^\Delta \stackrel{1}{\subseteq }{\frak{N}}^\Delta
\text{,}
\]
and we obtain a contradiction.$\blacksquare $\newline

%\bigskip

\begin{thebibliography}{99}
\bibitem{df}  \textsf{A. Defant and K. Floret}, \textit{Tensor norms and
operator ideals}, North - Holland Amsterdam, London, New York, Tokio 1993.

\bibitem{djt}  \textsf{J. Diestel, H. Jarchow, and A. Tonge}, \textit{%
Absolutely Summing Operators}, Cambridge University Press 1995.

\bibitem{gl}  \textsf{J. E. Gilbert and T. Leih}, \textit{Factorization,
tensor products and bilinear forms in Banach space theory}, Notes in Banach
spaces, pp. 182 - 305, Univ. of Texas Press, Austin, 1980.

\bibitem{glr}  \textsf{Y. Gordon, D. R. Lewis, and J. R. Retherford},
\textit{\ Banach ideals of operators with applications}, J. Funct. Analysis
14 (1973), 85 - 129.

\bibitem{gr}  \textsf{A. Grothendieck}, \textit{{R\'{e}sum\'{e} de la
th\'{e}orie m\'{e}trique des produits tensoriels topologiques}}, Bol. Soc.
Mat. S\~{a}o Paulo 8 (1956), 1 - 79.

\bibitem{h1}  \textsf{J. Harksen}, \textit{Tensornormtopologien},
Dissertation, Kiel 1979.

\bibitem{jo}  \textsf{H. Jarchow and R. Ott}, \textit{On trace ideals},
Math. Nachr. 108 (1982), 23 - 37.

\bibitem{l}  \textsf{H. P. Lotz}, \textit{Grothendieck ideals of operators
in Banach spaces}, Lecture notes, Univ. Illinois, Urbana, 1973.

\bibitem{oe1}  \textsf{F. Oertel}, \textit{Konjugierte Operatorenideale und
das }$\frak{A}$\textit{{-lokale Reflexivit\"{a}tsprinzip}}, Dissertation,
Kaiserslautern 1990.

\bibitem{oe2}  \textsf{F. Oertel}, \textit{Operator ideals and the principle
of local reflexivity}; Acta Universitatis Carolinae - Mathematica et Physica
33, No. 2 (1992), 115 - 120.

\bibitem{oe3}  \textsf{F. Oertel}, \textit{Composition of o{perator ideals
and their regular hulls}}; Acta Universitatis Carolinae - Mathematica et
Physica 36, No. 2 (1995), 69 - 72.

\bibitem{oe4}  \textsf{F. Oertel}, \textit{Local properties of accessible
injective operator ideals}; Czech. Math. Journal, 48 (123) (1998), 119-133.

\bibitem{oe5}  \textsf{F. Oertel}, \textit{Extension of finite rank
operators and local structures in operator ideals}; submitted

\bibitem{p1}  \textsf{A. Pietsch}, \textit{Operator ideals}, North - Holland
Amsterdam, London, New York, Tokio 1980.

\bibitem{pi}  \textsf{G. Pisier}, \textit{Factorization of linear operators
and geometry of Banach spaces}; CBMS Regional Conf. Series 60, Amer. Math.
Soc. 1986.
\end{thebibliography}
\end{document}